\documentclass[letterpaper,11pt]{article}

\usepackage{amsthm}
\usepackage{amsmath}
\allowdisplaybreaks[3]

\newtheorem{thm}{Theorem}[section]

\newtheorem{lemma}[thm]{Lemma}
\newtheorem{prop}[thm]{Proposition}

\newtheorem{conj}[thm]{Conjecture}

\newcommand{\beq}[1]{\begin{equation}\label{#1}}
\newcommand{\enq}[0]{\end{equation}}

\newcommand{\bn}[0]{\bigskip\noindent}
\newcommand{\mn}[0]{\medskip\noindent}
\newcommand{\nin}[0]{\noindent}

\newcommand{\sub}[0]{\subseteq}
\newcommand{\sm}[0]{\setminus}
\renewcommand{\dots}[0]{,\ldots,}

\newcommand{\A}[0]{{\cal A}}
\newcommand{\B}[0]{{\cal B}}

\newcommand{\m}[0]{{\cal M}}
\newcommand{\M}[0]{{\cal M}}

\newcommand{\ra}[0]{\rightarrow}

\newcommand{\II}[0]{{\bf I}}
\newcommand{\JJ}[0]{{\bf J}}

\newcommand{\XX}[0]{{\bf X}}

\newcommand{\E}[0]{{\sf E}}
\newcommand{\0}[0]{\emptyset}

\renewcommand{\qed}[0]{\begin{flushright} \rule{2mm}{3mm} \end{flushright}}

\newcommand{\C}[0]{\binom}
\newcommand{\Cc}[0]{\tbinom}
\newcommand{\ga}[0]{\alpha }

\newcommand{\gD}[0]{\Delta }
\newcommand{\gG}[0]{\Gamma }

\newcommand{\eps}[0]{\varepsilon }

\newcommand{\cov}[0]{~\cdot\hspace{-0.09in}>}

\newcommand{\sugg}[1]{}

\newcommand{\ma}[0]{{\rm ma}}
\newcommand{\mis}[0]{{\rm mis}}
\newcommand{\rr}[0]{r}
\newcommand{\ssss}[0]{s}

\usepackage[usenames]{color}

\begin{document}

\renewcommand{\thefootnote}{\fnsymbol{footnote}}
\footnotetext{AMS 2010 subject classification:  05C69, 05C35, 06A07}
\footnotetext{Key words and phrases:  maximal independent sets,
maximal antichains,
asymptotic enumeration, entropy, Shearer's Lemma}
\title{Counting maximal antichains and independent sets}
\author{L. Ilinca, J. Kahn}
\date{}

\maketitle

\begin{abstract}
Answering several questions of
Duffus, Frankl and R\"odl,
we give asymptotics for the logarithms of
(i)  the number of maximal antichains in the
$n$-dimensional Boolean algebra and
(ii) the numbers of maximal independent sets in the
covering graph of the $n$-dimensional hypercube
and certain natural subgraphs thereof.
The results in (ii) are implied by more general upper bounds on
the numbers of maximal independent sets in regular and biregular graphs.

We also mention some stronger possibilities involving
actual rather than logarithmic
asymptotics.

\end{abstract}

\section{Introduction}

Write ${\ma}(P)$ for the number of
maximal antichains of a poset $P$ and
$\mis(G)$ for the number of maximal independent sets in a graph $G$.
(For ``antichain" and ``independent set" see e.g. \cite{Engel} and
\cite{Diestel} respectively.)
Denote by $\B_n$ the Boolean algebra of order $n$
(that is, the collection of subsets of $\{1\dots n\}$ ordered by containment);
by $Q_n$ the ``covering" (or ``Hamming") graph of the $n$-cube
(the graph with vertex set
$\{0,1\}^n$ and two vertices adjacent if they differ
in exactly one coordinate); and by $\B_{n,k}$ the subgraph of
$Q_n$ consisting of strings of weight $k$ and $k+1$
(where, of course, {\em weight} means number of 1's).

We are interested in estimating the logarithms
of the quantities
$\ma (\B_n)$, ${\rm mis}(\B_{n,k}) $
and
${\rm mis}(Q_n) $, all problems
suggested by Duffus, Frankl and R\"odl in \cite{DFR1} and \cite{DFR2}.
As they observe
({\em cf.} the paragraph
following Conjecture \ref{conjdfr} below),
it is not hard to see that (with $\log =\log_2$ throughout)
\beq{lbs}
\log \ma (\B_n) \geq \Cc{n-1}{\lfloor n/2 \rfloor},
~
\log {\rm mis}(\B_{n,k}) \geq \Cc{n-1}{k}
~~
\mbox{and}~
\log {\rm mis}(Q_n) \geq 2^{n-2}.
\enq
On the other hand they show
\beq{kleit}
\log \ma (\B_n) < (1+o(1)) \Cc{n}{\lfloor n/2 \rfloor},
\enq
\[
\log {\rm mis}(Q_n) < (0.78+o(1))2^{n-1}
\]
and
\[
\log {\rm mis}(\B_{n,k}) < (1.3563+o(1))\Cc{n-1}{k}.
\]
Note that
$\log {\rm mis}(Q_n) \leq 2^{n-1}$
and
$\log {\rm mis}(\B_{n,k}) \leq \min\{\Cc{n}{k},\Cc{n}{k+1}\}$
are trivial
(see Proposition \ref{basicprop}(a)), while
\eqref{kleit} is Kleitman's celebrated
bound \cite{Kleitman} on the {\em total} number of antichains in $\B_n$.
In particular \eqref{kleit} makes no use of maximality,
and the authors of
\cite{DFR2} say:
``... the problem we are most interested in solving is this:
show there exists $\ga <1$ such that
$\log_2\ma (n) \leq \ga \C{n}{n/2}$"
(where their $\ma (n)$ is our $\ma(\B_n)$).

Here we settle all these problems, showing
that the lower bounds in \eqref{lbs}
are asymptotically tight, {\em viz.}
\begin{thm}
\label{TAnti}
{\rm (a)}
$\log \ma(\B_n) = (1+o(1))\tbinom{n-1}{\lfloor n/2\rfloor}$;

\mn
{\rm (b)}  $\log {\rm mis}(Q_n) = (1+o(1))2^{n-2}$;

\mn
{\rm (c)}  $ \log {\rm mis}(\B_{n,k}) = (1+o(1))\tbinom{n-1}{k}$
\end{thm}
\nin
(where $o(1)\ra 0$ as $n\ra \infty$).
Note that, with
$\B_{n,k}$ regarded as
the poset consisting of levels $k$ and $k+1$ of $\B_n$,
(c) also says
$\log {\rm ma}(\B_{n,k}) = (1+o(1))\tbinom{n-1}{k}$.

\medskip
The results giving (b) and (c) are actually far more general:
\begin{thm}
\label{mis_thm}
{\rm (a)}
For any $d$-regular, n-vertex graph $G$,
\begin{equation}
\log {\rm mis}(G) < \left\{\begin{array}{ll}
(1+o(1))\frac{n}{4}
&\mbox{if $G$ is triangle-free,}\label{mis_d_tr_free}\\
(1+o(1))\frac{n \log 3}{6} &\mbox{in general,}
\end{array}\right.
\enq
where $o(1)\ra 0$ as $d\rightarrow \infty$.

\medskip
\noindent
{\rm (b)} For any $(r,s)$-biregular, n-vertex (bipartite) graph $G$,
\[
\log {\rm mis}(G) < (1+o(1))\frac{rsn}{(r+s)^2},
\]
where $o(1)\ra 0$ as $\max\{r,s\}\ra\infty$.
\end{thm}

\nin
The first and third bounds are easily seen to imply (respectively)
parts (b) and (c) of Theorem \ref{TAnti}.
All three bounds 
are
best possible at the level of the asymptotics of the logarithm.
(This is shown, for example, by
\eqref{lbs}
for the first and third bounds, and by the graphs $H_n$ described in
Section \ref{Further}
for the second.)

\medskip
The rest of the paper is organized as follows.
Section \ref{Preliminaries} recalls what little background we need,
mainly a few known bounds on the parameter ``mis" plus
Shearer's entropy lemma.
The proofs of Theorems \ref{mis_thm}
and \ref{TAnti}(a) are given (in reverse order because the former is easier)
in Sections \ref{mis_proof} and \ref{PTa} respectively.
Finally, Section \ref{Further} suggests a
strengthening of Theorem \ref{TAnti}
and fills in
examples---possibly extremal---for
the second bound in Theorem \ref{mis_thm}(a).

The proofs of the theorems turn out not to involve too much
work once one gets on the right track.
In each case, seeking to identify an unknown (maximal independent set
or antichain) $I$,
we show that one can, at the cost of
specifying a few small subsets of $I$ and $\bar{I}$,
reduce determination of $I$ to determination of $I\cap Z$
for a relatively easily manageable subset $Z$ of our universe
(i.e. $V(G)$ or the ground set of $\B_n$).
Since the specified sets are small, there are not many ways to
choose them, so that possibilities for $I\cap Z$ contribute
the main term (plus part of the error) in each of our bounds.
Informally we tend to think of paying a small amount of ``information"
to reduce specification of $I$ to specification of $I\cap Z$.

Our original arguments for Theorems \ref{TAnti}
and \ref{mis_thm} were based on a
simple but (we think) novel combination of random sampling and entropy.
The proof of Theorem \ref{mis_thm} in Section
\ref{mis_proof} uses this approach.
As we recently noticed, there is an even simpler idea, due to
Sapozhenko \cite{Sap}, that can substitute for the
sampling/entropy part of the original argument; this
is explained at the end of
Section \ref{mis_proof}.
We have retained the original proof, in the case
of Theorem \ref{mis_thm} at least,
because we think the approach is interesting
and potentially useful elsewhere;
but in the case of Theorem \ref{TAnti},
for brevity's sake, we have omitted the original
and
given only the ``Sapozhenko version."

\section{Preliminaries}\label{Preliminaries}

The assertions of Theorem \ref{mis_thm} will reduce to the following
known and/or easy facts.

\begin{prop}\label{basicprop}

\mn
{\rm (a)}
For $G$ bipartite with bipartition $A\cup B$,
$$
{\rm mis}(G) \leq 2^{\min\{|A|,|B|\}}.
$$

\mn
{\rm (b)}  For any n-vertex graph G,
$$
{\rm mis}(G) \leq 3^{n/3},
$$
with equality iff $G$ is the disjoint union of $n/3$ triangles.

\mn
{\rm (c)}  For any n-vertex, triangle-free graph G,
$$
{\rm mis}(G) \leq 2^{n/2},
$$
with equality iff $G$ is a perfect matching.
\end{prop}
\nin
Here (a) is trivial; (b) was proved by
Moon and Moser \cite{MM} in answer to a question of Erd\H os and Moser;
and (c) is due to
Hujter and Tuza \cite{HT}.
We will not need the information about cases of equality.
Note that the $d$-regularity (with $d$ large) in Theorem \ref{mis_thm}(a)
multiplies the exponents in the
bounds directly implied by parts (b) and (c)
of Proposition \ref{basicprop} by roughly 1/2.

\medskip
We also need the following lemma of J. Shearer \cite{CFGS}.
(See e.g. \cite{McEliece} for entropy basics or \cite{Kahn}
for a quicker introduction and another application of Shearer's Lemma.)
We use $H$ for binary entropy and,
for a random vector $\XX=(\XX_1, \ldots, \XX_n)$ and $A\subset [n]$, set $\XX_A=(\XX_i:i\in A)$.
\begin{lemma} \label{Shearer}
Let $\XX=(\XX_1, \ldots, \XX_n)$ be a discrete random vector and $\mathcal{A}$ a collection of subsets (possibly with repeats) of $[n]$, with each element of $[n]$ contained in at least $m$ members of $\mathcal{A}$. Then
\[
H(\XX) \leq \frac{1}{m} \sum_{A\in \mathcal{A}} H(\XX_A).
\]
\end{lemma}
\nin
(The statement in \cite{CFGS} is less general,
but its proof gives Lemma \ref{Shearer}.)

\medskip
Finally, we write $\binom{n}{<t}$ for $\sum_{j<t} \binom{n}{j}$,
recalling that for $t\leq n/2$
(see e.g. \cite[Lemma 16.19]{Flum-Grohe})
\beq{binbd}
\Cc{n}{<t} \leq 2^{H(t/n)n}.
\enq

\section{Proof of Theorem \ref{mis_thm}} \label{mis_proof}

\mn
{\em Notation.}
For $G$ as in either part of the theorem,
we use $\gD=\gD(G)$ for maximum degree, ``$\sim$" for adjacency,
$N_x$ for the neighborhood of $x$, and $d(x)$ for $|N_x|$.
For
$S,T\subseteq V=V(G)$ and $x\in V$:
$d_S(x)=|N_x\cap S|$;
$\gG(S)=\left(\cup_{x\in S} N_x\right) \cup S$; $E(S)$ is the
sets of edges contained in $S$; and
$E(S,T)$ (used here only with $S\cap T\ =\0$) is
the set of edges having one end in each of
$S,T$.
We use $\m$ for the set of maximal independent sets of $G$.

\mn
{\em Proof.}
We work with a parameter $t$ to be specified below
and set $k = t^2\log \gD$.
For $I\in \M$, set $C_I=\{x\in V: d_I(x) \geq t\ln t\}$.
We may associate with each
$I\in \M$ some $I_0 \subseteq I$ of size $\lceil |I|/t\rceil$ with
\beq{I_0}
|C_I\sm \gG(I_0)|< n/t,
\enq
existence of such an $I_0$ being given by the observation that,
for $\II_0$ chosen uniformly from the $\lceil|I|/t\rceil$-subsets of $I$,
\begin{eqnarray*}
\E|C_I \sm \gG(\II_0)|&=&\sum_{v\in C_I} \Pr (v\not\in \Gamma(\II_0)) \\
&<& n (1-1/t)^{t\ln t}
< n/t.
\end{eqnarray*}

Let $X = V\sm (C_I\cup \gG(I_0))$, $Y=\{v\in X: d_X(v) \geq k\}$,
$J=I\cap Y$ and $Z= X\sm (Y\cup \gG(J))$.
We may choose $I\in \m
$ by specifying:
(i)  $I_0$;
(ii)  $C_I\sm \Gamma(I_0)$;
(iii)  $J$; and
(iv)  $I\cap Z$; so we just have to bound the numbers of
choices for these steps.  First, by
\eqref{binbd}, the number of choices in each of (i), (ii) is at most
$\exp_2[H(1/t)n]$.
Second, since $J\sub Y$ satisfies
\beq{Jsat}
|J\cap N_v|< t\ln t ~~\forall v\in X,
\enq
the log of the number of possibilities in (iii)
is at most $H(\JJ)$, where $\JJ$ is chosen uniformly from the collection
of subsets
$J$ of $Y$ satisfying \eqref{Jsat}.
Here Shearer's Lemma (with $\XX$ the indicator
of $\JJ$ and $\A=\{N_v\cap Y:v\in X\}$)
gives, using the definition of $Y$,
\begin{align}
H(\JJ)
&~\leq ~k^{-1}\sum_{v\in X} H(\JJ \cap N_v )
~\leq ~k^{-1}\sum_{v\in X} \log \Cc{\gD}{< t\ln t}
\nonumber\\
& ~<~ (k^{-1} t\ln t \log \gD)n
~=~(n\ln t)/t.
\label{X_1}
\end{align}

It remains to bound the number of possibilities in (iv).
Note that $I\cap Z$ is a {\em maximal} independent subset
of $Z$ (since $Z\cap \gG(I\sm Z)=\0$).
The (easy) point here is that the requirement
\beq{dZ}
d_Z(v)<k ~~\forall v\in Z
\enq
(implied by $Z\sub X\sm Y$)
limits the size of $Z$ (in (a))
or of the intersection of $Z$
with one of the parts of the bipartition
(in (b)).
We now consider these cases separately.

\mn
(a) From \eqref{dZ} we have
\[
d|Z| = \sum_{v\in Z} d(v) = 2|E(Z)|+ |E(Z, V\sm Z)| < k|Z|+d(n-|Z|),
\]
whence (note we will have $k<d$)
\[
|Z| < \frac{nd}{2d-k}
= (1/2 + O(k/d))n.
\]

\mn
Parts (c) and (b) of Proposition \ref{basicprop} thus
bound the log of the number of possibilities for $I\cap Z$ by
$(1+O(k/d))n/4$ if $G$ is
triangle-free, and by $(1+O(k/d))(\log 3)n/6$ in general.
Combining this with our bounds for (i)-(iii)
and setting $t=d^{1/3}$, we have
Theorem \ref{mis_thm}(a) with
$o(1) =O(\max\{(\log t)/t,(t^2\log d)/d\}) =O(d^{-1/3}\log d)$.

\mn
(b)
Let
the bipartition of $G$ be $A\cup B$, with $d(x)=r$ for $x\in A$ and,
w.l.o.g.,
($\gD=$) $r\geq s$.
Notice to begin that Proposition \ref{basicprop}(a) gives
\beq{bd1}
\log \mis(G)\leq |A| = sn/(r+s)  =(1+s/r)rsn/(r+s)^2.
\enq
We will prove the statement in (b) with
$o(1)=O(\min\{r^{1/3}s^{-2/3}\log r, s/r\})$,
so in view of \eqref{bd1}
may assume (to deal with a {\em very} minor detail below) that
\beq{pickiest}
r^{4/3}\log r < s^{5/3}.
\enq


Let $U=Z\cap A$ and $W=Z \cap B$.
We have
\[
r|U|=\sum_{v\in U} d(v) = \sum_{v\in U}d_{W}(v)+|E(U,B\sm W)|< k|U|+s(|B|-|W|),
\]
whence $r(1-k/r)|U|+s|W|< rsn/(r+s)$, implying either $|U|\leq (1-k/r)^{-1}rsn/(r+s)^2=(1+O(k/r))rsn/(r+s)^2$
or $|W|\leq rsn/(r+s)^2$
(where we used \eqref{pickiest} to say $k<r$).
Proposition
\ref{basicprop}(a) now bounds the number of possibilities for $I\cap Z$
by $\exp_2[(1+O(k/r))rsn/(r+s)^2]$.

Finally, combining this with \eqref{bd1} and
our earlier bounds for (i)-(iii), and setting $t=r^{2/3}s^{-1/3}$, gives
Theorem \ref{mis_thm}(b) with $o(1)$ on the order of
$$\min\left\{\max\left\{\frac{(\log t)(r+s)^2}{trs},\frac{k}{r}\right\},
\frac{s}{r}\right\}
=\Theta\left(\min\left\{\frac{r^{1/3}\log r}{s^{2/3}},
\frac{s}{r}\right\}\right).$$
\qed

\nin
{\em Sapozhenko's method}.
As promised, we next show how the first two paragraphs
of the above argument (through \eqref{X_1})
can be replaced by a remarkably
simple idea of
A. Sapozhenko
(see \cite[Theorem 6]{Sap} or, for another description,
\cite[Lemma 2.3]{Galvin}).

We again work with a parameter, say $b$, whose value will
be specified below.
For a given maximal independent set $I$,
let $X_1=V$ and repeat for $i=1,\ldots$ until no longer possible:
choose $x_i\in X_i\cap I$ with $d_{X_i}(x_i)\geq b$,
and set $X_{i+1}=X_i\sm (\{x_i\}\cup N_{x_i})$.
Let $Y=X_{q+1}$ be the final $X_i$, and notice that
$Y=V\sm \gG(\{x_1\dots x_q\})$ and
$$
d_Y(x)<b ~~\forall x\in Y\cap I.
$$
Set $Z =\{x\in Y:d_Y(x)< b\}$ ($\supseteq Y\cap I$).
We have

\mn
(i) $q<n/b$ (trivially);

\mn
(ii) $d_Z(x)<b  ~~\forall x\in Z$;  and

\mn
(iii) $I\cap Z$ is a maximal independent subset of $Z$
(since $I\cap (V\sm Z) = \{x_1\dots x_q\}\not\sim Z$).

\medskip
Now
(ii) corresponds to \eqref{dZ}, and the discussion under
(a) and (b) above (using (iii))
bounds the number of possibilities for $I\cap Z$
as before,
with the $k$'s replaced by $b$'s.
(For example, the main
bound in (b) is now $\exp_2[(1+O(b/r))rsn/(r+s)^2)]$.)

Essentially optimal values for $b$ are then
$(d\log d)^{1/2}$ and $rs^{-1/2}\sqrt{\log r}$
in (a) and (b)
respectively,
yielding $o(1)$'s
(as in the statement of the theorem) on the order of
$d^{-1/2}\sqrt{\log d}$ in (a), and, in (b),
$$\min\left\{\max\left\{\frac{(\log b)(r+s)^2}{brs},\frac{b}{r}\right\},
\frac{s}{r}\right\}
=\Theta\left(\min\left\{\sqrt{\frac{\log r}{s}},
\frac{s}{r}\right\}\right).$$
(So this also gives somewhat better error terms
than the original argument, though the errors are in any
case not likely to be close to the truth.)

\section{Proof of Theorem \ref{TAnti}(a)}\label{PTa}

\bn
This requires a little more care than the proof of Theorem \ref{mis_thm},
though the basic idea is similar.
As noted earlier, we skip our original argument and
just give the one based on ``Sapozhenko's method."

\mn
{\em Notation}. We write $\B$ for $\B_n$
(and follow a common abuse in using the same symbol for a
poset and its ground set).
Elements of $\B$ (usually denoted $x,y,z$)
may be thought of as either binary strings or subsets of $[n]$
(so for $x\in \B$ thought of as a string, $|x|$ is the number of 1's in $x$).
Set (for $i\in \{0\}\cup [n]$)
$L_i = \{x\in \B:|x| = i\}$, $\ell_i= |L_i|$ ($=\binom{n}{i}$) and,
for $S\sub \B$, $S_i=S\cap L_i$.
We also set $N=2^n$ and
$M=\binom{n}{\lfloor n/2\rfloor}$ ($=\Theta(n^{-1/2}N)$).

For $x\in \B$, $N_x$ is the neighborhood of $x$ in the
{\em comparability graph}, say $G$, of $\B$
(that is, the graph with $x\sim y$ iff $x<y$ or $x>y$).
Of course an antichain of $\B$ is just an independent set of $G$,
but in this case Theorem \ref{mis_thm} gives only a weak bound
on $\mis(G)=\ma(\B)$.
We write (as usual) $x\cov y$ ($x$ {\em covers} $y$) if
$x>y$ and there is no $z$ with $x>z>y$, and set
$d^+_S(y) =|\{x\in S: x\cov y\}|$ ($S\sub \B, y\in \B$).
Also for $S\sub \B$, we write
$\gG^+(S)$ for $\{y\in \B\sm S:\exists x\in S, x< y\}$
and $\ma(S)$ for the number of maximal antichains of
$S$ (more properly, of the restriction of $\B$ to $S$).
Finally, we
recall that $S\sub \B$
is {\em convex} if $x<y<z$ and $x,z\in S$ imply $y\in S$.

\mn
{\em Proof.}
Set $b= n^{3/4}\sqrt{\log n}$.
Given a maximal antichain $I$,
let $X_1=\B$ and repeat for $i=1,\ldots$ until no longer possible:
choose $x_i\in I\cap X_i$ with $|N_{x_i}\cap X_i|\geq b$,
and set $X_{i+1}=X_i\sm (\{x_i\}\cup N_{x_i})$.
Let $Y=X_{q+1}$ be the final $X_i$---so in particular
$
d^+_Y(y)< b ~~\forall y\in I\cap Y
$---and set $Z =\{y\in Y:d^+_Y(y)< b\}$ ($\supseteq I\cap Y$).

The number of possibilities for $\{x_1\dots x_q\}$
is at most $2^{H(1/b)N}$ (since
$q<N/b$), and we have
\beq{d+z}
d^+_Z(x)< b  ~~\forall x\in Z,
\enq
$$
\mbox{{\em $I\cap Z$ is a maximal antichain of $Z$}}
$$
(since $I\cap (\B\sm Z) = \{x_1\dots x_q\}\not\sim Z$),
and, we assert,
\beq{cvx'}
\mbox{$Z$ is convex.}
\enq

\mn
{\em Proof of} \eqref{cvx'}.
Since $Y= \B\sm \bigcup_i(\{x_i\}\cup N_{x_i})$ is obviously convex,
\eqref{cvx'} follows from the observation that
$$\mbox{$Y\sm Z$ is a downward-closed subset of $Y$}.$$
For suppose that---now regarding elements of $\B$ as
subsets of
$[n]$, for which we prefer capitals---$A\in Y$ and $A\sub B\in Y\sm Z$.
Since $B\not\in Z$, there are
distinct $i_1\dots i_{b}\in [n]\sm B$
with $B\cup \{i_j\}\in Y$
$\forall j\in [b]$,
whence, since $Y$ is convex (and $A\in Y$),
$A\cup \{i_j\}\in Y$ $\forall j\in [b]$.
But then $A\in Y\sm Z$.\qed

Theorem \ref{TAnti}(a) (with $o(1)$ on the order of
$n^{-1/4}\sqrt{\log n}$) thus follows from the next two assertions
(and the fact that $M\leq 2\C{n-1}{\lfloor n/2\rfloor}$).

\mn
{\em Claim} 1.
If $Z\sub \B$ is convex and satisfies \eqref{d+z},
then
\beq{ZM}
|Z| < (1+O(b/n))M.
\enq

\mn
{\em Claim} 2.
If $Z\sub \B$ is convex, then
$\ma(Z)\leq 2^{|Z|/2}$.

\mn
(Note this does require convexity; e.g. it fails if $Z$
is a 3-element chain.)

\mn
{\em Proof of Claim} 1.
We may assume $Z\sub L_{\rr}\cup\cdots \cup L_{\ssss}$, with
$(\rr,\ssss)=(.4 n,.6n)$
(since the rest of $\B$ is too small to affect \eqref{ZM}).
Let
$F=\gG^+(Z)$,
$f_i=|F_i|$ and $z_i=|Z_i|$.
The degree assumption \eqref{d+z}
implies that, for any $i\in \{\rr\dots \ssss\}$,
\[
(i+1)f_{i+1}\geq \partial (F_i\cup Z_i,F_{i+1})\geq
(n-i)f_i + (n-i-b)z_i
\]
(where $\partial(S,T):=|\{(x,y):x\in S,y\in T,x<y\}|$),
or, since $(n-i)/(i+1) = \ell_{i+1}/\ell_i$,
\beq{fili}
f_{i+1}\geq
\frac{\ell_{i+1}}{\ell_i}f_i +
\left(\frac{\ell_{i+1}}{\ell_i}-\frac{b}{i+1}\right)z_i
\geq \frac{\ell_{i+1}}{\ell_i}[f_i+(1-\eps)z_i],
\enq
with $\eps = 2.5b/n$.
Composing the inequalities \eqref{fili} (for $i=\rr\dots \ssss$) gives
\[
f_{\ssss+1}\geq (1-\eps)\sum_{j=\rr}^{\ssss}\frac{\ell_{\ssss+1}}{\ell_j}z_j,
\]
so that
\[
\frac{|Z|}{M}~\leq ~ \sum_{j=\rr}^{\ssss}\frac{z_j}{\ell_j}~\leq ~
(1-\eps)^{-1}\frac{f_{\ssss+1}}{\ell_{\ssss+1}}~\leq ~(1-\eps)^{-1}.
\]
\qed

\nin
{\em Proof of Claim} 2.
This is an induction along the lines of (but easier than)
the argument of \cite{HT}.
Set $|Z|=m$.
If $Z$ does not contain a chain of length 3, then the
comparability graph of $Z$ is bipartite
and we may apply Proposition \ref{basicprop}(a).
Otherwise there are $x,y\in Z$ with $x<y$ and $|y|\geq |x|+2$,
whence, since $Z$ is convex, $|\gG^+(x)\cap Z|\geq 3$.
We may, of course, further require that $x$ be minimal
in $Z$ (so $Z\sm \{x\} $ is convex), and then
induction gives
$$\ma(Z)\leq \ma(Z\sm \{x\}) + \ma(Z\sm (\{x\}\cup \gG^+(x)))
\leq 2^{(m-1)/2} + 2^{(m-4)/2} < 2^{m/2}.$$
(Notice that $\ma(Z\sm (\{x\}\cup \gG^+(x)))$ bounds
the number of maximal antichains of $Z$ containing $x$, since
for each such $A$, $A\sm \{x\}$ is contained in some maximal
antichain $A'$ of $Z\sm \{x\}$, and is recoverable from $A'$
{\em via} $A=(A'\cup \{x\})\sm \{y\in Z:y\sim x\}$.)\qed

\section{A stronger conjecture}\label{Further}

In closing we would like to suggest that it might be
possible to give the actual asymptotics
(rather than just the asymptotics of the log) for the quantities
considered in Theorem \ref{TAnti}.  This does not look
easy, but we can at least guess what the truth should be:

\begin{conj}\label{conjdfr}
{\rm (a)}
$ \displaystyle
\ma(\B_n)=\left\{\begin{array}{ll}
(1+o(1))\,n \,\exp_2[\binom{n-1}{(n-1)/2}] &\mbox{if $n$ is odd,}\\
(2+o(1))\,n \,\exp_2[\binom{n-1}{n/2}] &\mbox{if $n$ is even};
\end{array}\right.
$

\mn
{\rm (b)}
$ \mis (Q_n)=(2+o(1))\,n \,\exp_2[2^{n-2}]$;

\mn
{\rm (c)}
$ \mis(\B_{n,k})=(1+o(1))\,n \,\exp_2[\binom{n-1}{k}]$

\mn
(where $o(1)\ra 0$ as $n\ra \infty$).
\end{conj}
\nin
The easy lower bounds
are based on the observation (from
\cite{DFR1})
that for any graph $G$
and induced matching $M$ of $G$,
each of the
$
2^{|M|}
$
sets consisting of one vertex from each edge of $M$
extends to at least one maximal independent set, and these
extensions are all different.
For example, the lower bound for $n$ odd
in Conjecture \ref{conjdfr}(a)
(other cases are similar)
is obtained by noting that,
for each $i\in[n]$,
the set of pairs
$$
\{\{x,x^i\}:x\in \{0,1\}^n, x_i=0, |x|=(n-1)/2\}
$$
(where $x^i$ is gotten by flipping the $i$th coordinate of $x$)
is an induced matching in the comparability graph of $\B_n$,
and that there is only an
insignificant amount of repetition in the corresponding list of at
least $n\exp_2[\binom{n-1}{(n-1)/2}]$ maximal independent
sets.

\medskip
Finally  we give the promised construction for the second
bound in Theorem \ref{mis_thm}(a).
For
$d\equiv 2\pmod{3}$, let $T$ be the disjoint union of $(d-2)/3$ triangles; let
$H$ consist of two disjoint copies of $T$ plus all edges between them; and, for
$n$ divisible by $2(d-2) $, let $H_n$ be the union of $\frac{n}{2(d-2)}$
disjoint copies of $H$.  Then $H_n$ is $d$-regular with
$
\mis (H_n) =
2^{n/(2(d-2))}3^{n/6},
$
and it seems not impossible
that this is extremal:
\begin{conj}\label{Hconj}
For any d-regular, n-vertex G, $\mis(G)\leq 2^{n/(2(d-2))}3^{n/6}$.
\end{conj}
\nin
This would be analogous to the fact---recently proved in spectacularly
simple fashion by Yufei Zhao \cite{Zhao1} (but by reducing to
the bipartite case proved a little less simply in
\cite{Kahn})---that the {\em total} number of independent sets in such a $G$ is
at most $(2^{d+1}-1)^{n/(2d)} $,
a value achieved by a disjoint union of $K_{dd}$'s
whenever $2d|n$.
One would, of course, hope to also have analogues for the other parts of
Theorem \ref{mis_thm};
but we don't see good candidates for these,
and suspect that they do not have
clean answers.

\bigskip
\begin{flushleft}
Department of Mathematics\\
Indiana University\\
Bloomington IN 47405 USA\\
ilinca@indiana.edu\\
$~$\\
Department of Mathematics\\
Rutgers University\\
Piscataway NJ 08854 USA\\
jkahn@math.rutgers.edu
\end{flushleft}

\end{document}